\documentclass[a4paper,11pt]{article} 
\usepackage[utf8]{inputenc} 
\usepackage{lmodern}
\usepackage[T1]{fontenc}
\usepackage[english]{babel}
\usepackage{amsmath}
\usepackage{amsthm}
\usepackage{amsfonts}
\usepackage{makeidx}
\usepackage[nottoc, notlof, notlot]{tocbibind}
\usepackage{amssymb}
\usepackage{enumerate}
\usepackage{euscript}
\usepackage{lmodern}
\usepackage{parskip}
\usepackage[final]{showkeys} 

\newcommand*\diff{\mathop{}\!\mathrm{d}}

\theoremstyle{plain}
\newtheorem{theorem}{Theorem}[]

\newtheorem{lemma}[theorem]{Lemma}
\newtheorem{corollary}[theorem]{Corollary}
\newtheoremstyle{sty}
{9pt} 
{3pt} 
{} 
{} 
{\bfseries} 
{.} 
{.5em} 
{} 

\theoremstyle{sty}

\newtheorem{remark}{Remark}

\title{Expansion of an error threshold for a finite population in the Moran model}
\author{Maxime Berger
\footnote{
\noindent
D\'epartement de math\'ematiques et applications, \'Ecole normale sup\'erieure, CNRS, PSL University, 75005 Paris, France
}}
\date{\today}
\begin{document}
\maketitle

\begin{abstract}
\noindent
We propose a new definition for the error threshold of a population evolving through mutation and selection.
We compute the correction term due to the finiteness of the population by estimating the lifetime of master sequences.
Our technique consists in bounding from above and below the number of master sequences in the Moran model, by birth and death chains. The expectation of this lifetime is then computed with the help of explicit formulas which are in turn expanded with Laplace method.
The first term after~$\ln \sigma / \ell$ is computed, it scales as~$1/\sqrt{\ell m}$, where~$\ell$ is the genome size and~$m$ the number of individuals in the population.
\end{abstract}
 \noindent
\textit{keywords}: Markov chain, phase transition, genetics\\
 \noindent
\textit{Mathematic Subject Classification}: 60J10, 92D10

\section{Introduction}
Evolution is a macroscopic phenomenon that relies on two microscopic forces: mutation and selection.
In the 70s, Manfred Eigen \cite{eigenSelforganizationMatterEvolution1971} introduced a mathematical model to understand the evolution of a population of prebiotic molecules. He studied the structure of the population after a long time and proved the existence of a phase transition. If the parameter governing the mutation rate is above a certain error threshold~$q^*$, the population becomes completely random and all the genetic information is lost, if it is belown~$q^*$, a positive concentration of the population retains the fittest genotype. 
Let~$\ell$ be the length of the genome of the macromolecules, and~$q$ the probability for mutation per site, so that the product~$\ell q$ represents the average number of mutations for each reproduction event.
Eigen computed the error threshold~$q^*$ in a model with an infinite number of macromolecules, in an asymptotic regime where 
\[ \ell \rightarrow \infty\,, \quad q\rightarrow 0\,,\]
and he obtained
\[ q^* \mathop{\sim} \frac{\ln \sigma}{\ell}\,,\]
where~$\sigma$ stands for the fitness of the fittest genotype, all the others having fitness~$1$. This fitness landscape is known as the sharp-peak landscape,
it is commonly used in population genetics works, mainly because computations are easier, but it is also a plausible framework: in real life, most mutations do not modify the fitness.
The present work is also done within this framework.

However, real populations are not infinite and it is necessary to study models with a finite population.
Our goal in this text is to compute the correction term in the expansion of the error threshold due to the finite population size.
We study the Moran model \cite{moranRandomProcessesGenetics1958}, introduced in the 50s, which has been shown to converge to Eigen's model by Dalmau in \cite{dalmauConvergenceMoranModel2017}.
This convergence allows us to derive results on the quasispecies model working with population genetics models. As Wilke argued in \cite{wilkeQuasispeciesTheoryContext2005a}, this is a sensible strategy and it is supported by numerous simulations. Understanding the scaling with the size of the population in the quasispecies model leads to many applications like optimizing genetic algorithms~\cite{ochoaErrorThresholdsTheir1999} or finding ways to eradicate a population of viruses by increasing their mutation rate~\cite{tripathiStochasticSimulationsSuggest2012}.
Let us first describe the model.

We consider a population of~$m$ individuals whose genetic material is coded with a string of~$\ell$ characters chosen in~$\{1,\cdots,\, \kappa\}$. All of the~$\kappa^\ell$ genotypes have fitness equal to~$1$ except one sequence, say~$1\cdots 1$, which has fitness equal to~$\sigma$, with~$\sigma >1$. The sequence~$1\cdots 1$ is called the master sequence. 
At each generation, one individual is chosen to be a parent, but master sequences have a selective advantage: they 
are~$\sigma$ times more likely to be chosen, thus they have better chances to leave more offsprings.
All the other sequences are equally likely to be chosen. 
The chosen individual is replicated, yet the replication process is error prone, due to mutations, each bit of its genome is changed independently with probability~$q$ into one of the other~$\kappa-1$ letters.
The offspring replaces a randomly chosen individual in the population.
In particular, the size of the population is constant equal to~$m$.
This process is repeated indefinitely.

The Moran model, like the Eigen model, presents a phase transition separating a regime of chaos and a regime where master sequences occupy a non negligible proportion of the population.
However, the mathematical definition of the critical parameter is delicate and several choices are possible.
Suppose the population starts with no master sequences. Some random mutations happen and new genotypes are discovered without changing the global fitness, and the selection plays no role since every sequence has the same fitness.
This phase is called the \textbf{neutral phase}, the population evolves as if master sequences did not exist. 
At some point, a random mutation will discover a master sequence. From this moment on, the selection enters the game and it helps to keep the master sequences in the population.
This phase is called the \textbf{quasispecies phase}, the master sequences occupy a significant proportion of the population, along with a cloud of mutants consisitng of individuals that are genetically closed to the master sequence.
Increasing the mutation rate reduces the stability of the quasispecies phase because offsprings of master sequences are then less likely to be master sequences, however the stability of the neutral phase is barely changed.
The mutation probability can be chosen such that both phases are equally stable in some sense, this critical value is a first possible definition for the error threshold. 
It is also possible to increase this probability even more in such a way that the quasispecies phase will stop being stable at all, master sequences will then disappear very quickly. This point gives another possible definition of the critical parameter, analogous to the so called spinodal point in statistical mechanics.
In \cite{nowakErrorThresholdsReplication1989} Nowak and Schuster looked for a critical parameter in some modified version of the Moran model. They computed the stationary measure of the number of master sequences, this measure has two maxima: one corresponding to the neutral phase and another one to the quasispecies phase. When the mutation parameter is increased, the second maximum becomes smaller until it vanishes and only one maximum remains.
This point corresponds to our second definition, the spinodal point.
Nowak and Schuster computed the expansion of this point due to the finite population size, they obtained the following correction for~$q^*$:
\[ q^* \mathop{\sim} \frac{\ln \sigma}{\ell}
-\frac{2\sqrt{\sigma-1}}{\ell\sqrt{m}}\,.\]
However, they made an assumption that lead them to underestimate the number of master sequences, so the true critical point is larger in the Moran model.

We aim at a different goal here: to find the parameter at which both phases are equally stable.
For this task, we need to estimate the relative stability of both phases, in fact we will estimate the time needed to get from one phase to the other.

When in the neutral phase, the time~$\tau^*$ needed to discover a master sequence is called the \textbf{discovery time}. When in the quasispecies regime, the master sequences are present in the population during a certain time~$\tau_0$ that we call the \textbf{persistence time}.
It is by estimating these two times and comparing them that we will estimate our error threshold.

The expectation of the time~$\tau^*$ has been estimated in \cite{cerfCriticalPopulationError2015b}, it has been shown to be independent of the parameter~$q$ in first approximation and of order~$\kappa^\ell$, which is the size of the entire sequence space.
More precisely, in a regime where
\[
\ell \rightarrow \infty, \quad q\rightarrow 0, \quad m\rightarrow \infty\,,
\]
and such that
\[
\ell q\rightarrow a \in ]0, \infty[, \quad \frac{m}{\ell}\rightarrow \alpha\in [0, \infty]\, ,
\]
the expectation of the time~$\tau^*$ admits the following expansion:
\[
\lim \frac{1}{\ell}\ln E(\tau^*)= \ln \kappa\,.
\]
Our main goal here is to estimate the time~$\tau_0$. As the last sentence of \cite{eigenMolecularQuasiSpecies2007} states, estimates on the lifetime of the metastable quasispecies is crucial to understand the Eigen model. We need~$\tau_0$ with sufficient precision to decide whether the ratio~$E(\tau^*)/E(\tau_0)$ tends towards~$0$ or towards infinity. 
If the ratio~$E(\tau^*)/E(\tau_0)$ goes to zero, master sequences will occupy a positive proportion of the population, whereas if the ratio goes to infinity, master sequences will vanish. 

We always work in the following asymptotic regime:
\begin{equation}
\label{regime}
m\rightarrow \infty\,, \quad \ell \rightarrow \infty\,, \quad q\rightarrow 0\,,
\end{equation}
we must however add some conditions. First we assume that the product~$\ell q$ admits a finite and positive limit:
\begin{equation}
\label{regimeconditiona}
\ell q\rightarrow a\in ]0,\infty[\,. 
\end{equation}
Moreover, we must ensure that this limit~$a$ is not too large, because we will need the fact that
\begin{equation}
\label{regimeconditionrho}
\sigma e^{-a}-1 \geq 0\,\,.
\end{equation}
This quantity is related to the average number of master sequences in the quasispecies phase. The opposite inequality is also interesting. In fact, if~$\sigma e^{-a} < 1$, master sequences disappear much quicker and their lifetime is then much smaller. The way we handled our remainders does not allow us to attain the precision we would need.

We make a third hypothesis. We do not yet understand its relevance, but our computations become easier if we assume that
\begin{equation}
\label{regimeconditionsigma}
\sigma e^{-a} -1 \neq \sigma -2\,.
\end{equation}
The case when these two quantities are equal is dicussed in remark~\ref{rem}.

The estimation on the expectation of the persistence time is stated in theorem~\ref{theoremTemps}, it does require very few hypothesis on the asymptotic behaviour of the parameters~$\ell$,~$m$ and~$q$. The asymptotic expansion of the error threshold, however, depends largely on the choice of the relative sizes of parameter~$\ell$ and~$m$. 
Our main goal is to go further in the expansion~$ q^* = \frac{\ln \sigma}{\ell}$, which corresponds to the case~$a=\ln \sigma$, that is why we assume this equality  in the following theorem.

\vspace{10pt}
\begin{theorem}
\label{seuilErreur}
We suppose that $a=\ln \sigma$.
In the asymptotic regime~\eqref{regime} with conditions~\eqref{regimeconditiona},~\eqref{regimeconditionrho} and~\eqref{regimeconditionsigma},
if 
\[
\frac{m}{\ell}\rightarrow \infty,
\quad
\frac{\ell^2}{m\ln m}\rightarrow \infty\,,
\]
then the error threshold defined by the equality of the two times~$\tau_0$ and~$\tau^*$ expands as 
\[
q^* =\frac{\ln\sigma }{\ell} -\frac{\sqrt{2(\sigma-1)\ln \kappa} }{\sqrt{\ell m}}\,,
\]
in the following sense: if we take
\[q =\frac{\ln\sigma }{\ell} -\frac{c }{\sqrt{\ell m}},\]
with 

$\bullet$~$c < \sqrt{2(\sigma-1)\ln \kappa}$, then no master sequences are present in the population, 

$\bullet$~$c>\sqrt{2(\sigma-1)\ln \kappa}$, the population has a positive concentration of master sequences.
\end{theorem}

In the last case, the proportion of master sequences in the quasispecies phase is equivalent to
$$
\frac{c}{\sigma-1}\,\sqrt{\frac{\ell}{m}}\,.
$$
Contrary to Nowak and Schuster, our correction term scales as the inverse of the square root of the chain length~$\ell$, which leads indeed to a smaller error threshold.
The reason for this discrepancy is, as we explained earlier, that Nowak and Schuster compute the spinodal point of their model, which is higher than the critical threshold we consider here.

The previous expansion is derived from the expansion of the persistence time stated in the next theorem.

\vspace{10pt}
\begin{theorem}
\label{theoremTemps}
In the asymptotic regime~\eqref{regime} with conditions~\eqref{regimeconditiona},~\eqref{regimeconditionrho} and~\eqref{regimeconditionsigma},
the persistence time is of order
\[
E\Big(\tau_0 \, | \, N_0=1\Big)=
 \exp \Bigg( 
m\varphi\Big((1-q)^\ell\Big)
+ O\Big((1+mq)\ln m\Big)\Bigg)\,,
\]
where
\[
\varphi(x)= \frac{\sigma\big(1-x)\ln\frac{\sigma(1-x)}{\sigma -1}+\ln (\sigma x)}{1-\sigma(1-x)}\,.
\]
\end{theorem}

From the estimate on the time~$\tau_0$, we can also decide which of the two regimes is dominant according to some other choice of the parameters. Let us consider a regime where the size of the population is proportional to the genome length.
\vspace{10pt}
\begin{corollary}
\label{cor}
In the asymptotic regime~\eqref{regime}, with conditions~\eqref{regimeconditiona},~\eqref{regimeconditionrho} and~\eqref{regimeconditionsigma},
if 
\[\frac{m}{\ell} \rightarrow
\alpha\in]0,\infty[\,,
\]
then no error threshold exists if~$\alpha<{\ln \kappa}/{\ln\sigma}$, in this case the population is always neutral.

If~$\alpha>{\ln \kappa}/{\ln\sigma}$, there is an error threshold at the point~$q^*$ satisfying
\[
(1-q^*)^\ell = \varphi^{-1}\Big(\frac{\ln\kappa}{\alpha}\Big)\,.
\]

\end{corollary}
Since~$\sigma(1-q)^\ell-1$ is related to the average number of master sequences, the persistence time will be much smaller if this quantity is negative, 
so the concentration of the master sequences will still be negligible in this case.

The theorem~\ref{theoremTemps} on the time~$\tau_0$ permits us to describe the situation in a regime where the size of the population is much smaller than the genome length. 
In the asymptotic regime~\eqref{regime}, with conditions~\eqref{regimeconditiona},~\eqref{regimeconditionrho} and~\eqref{regimeconditionsigma},
in the case \[
\frac{m}{\ell}\rightarrow 0\,,
\]
there can be no error threshold, the regime is always neutral. The reason is that when there are fewer individuals, it is easier to randomly loose master sequences.
This is also the case when~$\sigma(1-q)^\ell-1$ is negative.

\begin{remark}
\label{rem}
In the asymptotic regime~\eqref{regime}, with conditions~\eqref{regimeconditiona}, and~\eqref{regimeconditionrho} only, and if we set
\[
\sigma e^{-a} -1 = \sigma -2\, ,
\]
the expectation of the persistence time~$\tau_0$ is estimated by
\begin{multline*}
E\Big(\tau_0 \, | \, N_0=1\Big)=
 \exp \Bigg( 
m\varphi\Big((1-q)^\ell\Big)\\
+ O\bigg(\frac{mq}{
L_q\big(\sigma L_q + (\sigma -1)q\big)}\bigg)
+ O\Big((1+mq)\ln m\Big)\Bigg)\,,
\end{multline*}
where
\[
\varphi(x)= \frac{\sigma\big(1-x)\ln\frac{\sigma(1-x)}{\sigma -1}+\ln (\sigma x)}{1-\sigma(1-x)}\,,
\]
and 
\[
L_q = \sigma - 2 - \Big(\sigma (1-q)^\ell -1\Big)\,.
\]
\end{remark}

Even if we concluded with the absence of error threshold in some regimes, other definitions of critical parameters are possible in these regimes. We could for example look at the probability~$q$ such that the time~$\tau_0$ grows smaller than a polynomial function of~$m$, or even the probability which would lead to keep~$\tau_0$ bounded as~$m$ increases.

We will perform the computations for two simpler models: one underestimating the number of master sequences, and another one overestimating it.
We obtain two bounds on the persistence time which yield two bounds on the error threshold.
\section{The birth and death process}

We follow the strategy of~\cite{nowakErrorThresholdsReplication1989} to simplify the original process, namely, we classify the individuals in only two types.
The first type~$T_1$ gathers all the master sequences, all other sequences are put in the second type~$T_2$.
Any master sequence undergoing mutation is put in the type~$T_1$.
In order to bound from above and below the error threshold, we modify the probability for a individual of type~$T_2$ to become a master sequence.
Since we study a regime where there are few mutations ($q\rightarrow 0$), it is easier to get master sequences if every individual needs only one mutation to become a master sequence.
Under this scheme, the number of master sequences will be greater than in the original process and thus we will obtain a longer persistence time.
On the contrary, if~$\ell$ mutations are needed, there will be less master sequences, and master sequences will extinct faster.
The two cases will be coded by the letter~$\theta$, which will represent the number of mutations needed to become a master sequence.
This notation will allow us to perform the associated computations only once: the case~$\theta =1$ will give an upper bound on the time~$\tau_0$, and the case~$\theta=\ell$ will give a lower bound.

The number of master sequences at time~$t$, denoted by~$N_t$, evolves according to a birth and death process, indeed:

$\bullet$ The number~$N_t$ increases by~$1$ when a master sequence does not mutate and replaces an individual of type~$T_2$ or when an individual of type~$T_2$ undergoes the right mutations to become a master sequence and replaces an individual of type~$T_2$. 
For~$k$ between~$0$ and~$m-1$, we denote by~$\delta_k$ the probability that~$N_t$ jumps from~$k$ to~$k+1$:
\begin{align}
\label{defDeltaprob}
&\forall t \geq 0\quad \forall k \in \{0,\cdots, m-1\} \nonumber\\
&\qquad\qquad\quad\delta_k
=\,P\Big(N_{t+1}=k+1\, \Big|\, N_t=k\Big) \nonumber\\
&\qquad\qquad\qquad=\,\frac{\sigma \frac{k}{m}\big(1-\frac{k}{m}\big) (1-q)^\ell+\big(1-\frac{k}{m}\big)^2(1-q)^{\ell-\theta}\big(\frac{q}{\kappa-1}\big)^\theta}{\sigma \frac{k}{m} +1 -\frac{k}{m}}\,.
\end{align}

$\bullet$ The number~$N_t$ decreases by~$1$ if a master sequence mutates and replaces a master sequence or if an individual of type~$T_2$ does not become a master sequence and replaces a master sequence. For~$k$ between~$1$ and~$m$, we denote by~$\gamma_k$ the probability that~$N_t$ jumps from~$k$ to~$k-1$:
\begin{align}
\label{defGammaprob}
&\forall t \geq 0\, \,\,\forall k \in \{1,\cdots, m\}\,\nonumber \\
&\qquad\quad\quad\gamma_k=
\,P\Big(N_{t+1}=k-1\, \Big|\,N_t=k\Big)\nonumber\\
&\qquad\qquad\quad=\, \frac{\sigma \big(\frac{k}{m}\big)^2 \Big(1-(1-q)^\ell\Big)+\frac{k}{m}\big(1-\frac{k}{m}\big)\Big(1-(1-q)^{\ell-\theta}\big(\frac{q}{\kappa-1}\big)^\theta\Big)}{\sigma \frac{k}{m}+1 -\frac{k}{m}}\,.
\end{align}

We define the persistence time~$\tau_0$ by the lifetime of master sequences:
\[
\tau_0=\inf \, \{\, t\geq0: N_t=0\, \}\,.
\]
For such birth and death processes, there exists an explicit formula for the expectation of the persistence time~$\tau_0$ stated in the next lemma.
We set~$\pi_0=1$ and
\[
\forall i \in \{1,\cdots, m\}\qquad \pi_i=\frac{\delta_1\cdot\cdot\cdot\delta_i}{\gamma_1\cdot\cdot\cdot\gamma_i}\,.
\]
In the sequel, we will use~$\pi_m / \delta_m$, which is defined even though~$\delta_m$ is not, indeed, we have
\[
\frac{\pi_m}{\delta_m} = \frac{\delta_1\cdot\cdot\cdot\delta_{m-1}}{\gamma_1\cdot\cdot\cdot\gamma_m}\,.
\]
\vspace{10pt}
\begin{lemma}\label{vieMort}
The expectation of the persistence time~$\tau_0$ started from~$N_0=1$ is given by
\[ 
E\Big(\tau_0 \, | \, N_0=1\Big)=\sum^m_{i=1}\frac{1}{\delta_i} \, \pi_{i}\,.
\]
\end{lemma}
We shall rely on this formula to compute the expectation of the persistence time. We will start by estimating~$\ln (\delta_k/\gamma_k)$, we will then focus on
\[
\ln \pi_i=\sum_{k=1}^{i}\,\ln \frac{\delta_k}{\gamma_k},\quad i\in\{1,\cdots,m\}\,,
\] 
add up the quantities~$\exp(\ln \pi_i)/\delta_i$, and implement Laplace's method to estimate the sum.

\subsection{Computation of~$\ln \delta_k / \gamma_k$}
In order to compute the ratio~$\frac{\delta_k}{\gamma_k}$, let us first factorise the expressions.
Let us introduce the notation
\begin{equation}
\label{Q}
Q = \frac{ q }{ (1-q)(\kappa-1) }\, .
\end{equation}
By factorising the probability~$\delta_k$ from expression~\eqref{defDeltaprob}, we obtain
\begin{equation}
\label{delta}
\delta_k = \Big(1-\frac{k}{m}\Big)\sigma (1-q)^\ell \frac{\frac{Q^\theta}{\sigma}+\big(1-\frac{Q^\theta}{\sigma}\big)\frac{k}{m}}{\sigma \frac{k}{m} +1 -\frac{k}{m}}\,.
\end{equation}
Similarly, we can rewrite~$\gamma_k$ thanks to expression~\eqref{defGammaprob} as
\[
\gamma_k = \frac{k}{m}\frac{ 1 - (1-q)^\ell Q^{\theta} +\Big( \sigma - 1 - \sigma (1-q)^\ell + (1-q)^\ell Q^{\theta} \Big)\frac{k}{m}}{\sigma \frac{k}{m}+1 -\frac{k}{m}}\,.
\]
We set~$L_q$ to be
\begin{equation}
\label{Lq}
L_q = \sigma - 1 - \sigma (1-q)^\ell\,,
\end{equation}
and we define two functions~$\psi$ and~$\phi$ by
\begin{equation}
\label{psi}
\psi(x) = \frac{ Q^{\theta} }{ \sigma }+\Big( 1 - \frac{ Q^{\theta} }{ \sigma } \Big) x \, ,
\end{equation}
and 
\begin{equation}
\label{phi}
\phi(x) = 1 - (1-q)^\ell Q^{\theta} +\Big( L_q + (1-q)^\ell Q^{\theta} \Big)x\, .
\end{equation}

With these notations, the ratio~$\delta_k/\gamma_k$ can be rewritten as
\[
\frac{\delta_k}{\gamma_k} = \,
\frac{1- \frac{k}{m}}{\frac{k}{m}}\, \sigma (1-q)^\ell\,
 \frac{\psi\big(\frac{k}{m}\big)}{\phi\big(\frac{k}{m}\big)}\,.
\]

Sums are easier to work with than products, so we start by looking for an estimate of
\begin{equation}
\label{lnd/g}
\ln \frac{\delta_k}{\gamma_k} = 
\ln\bigg( \frac{1-\frac{k}{m}}{\frac{k}{m}} \bigg) + \ln \Big(\sigma (1-q)^\ell\Big)+
\ln\psi\Big(\frac{k}{m}\Big) -
\ln\phi\Big(\frac{k}{m}\Big)\,.
\end{equation}
Let~$i$ be an integer in~$\{1,\cdots, m\}$, summing identity~\eqref{lnd/g} between~$1$ and~$i$ gives
\begin{align*}
\ln \pi_i = &
\,\sum_{k=1}^i \ln \frac{\delta_k}{\gamma_k}\\
=&
\,\sum_{k=1}^i \ln \left( \frac{ 1 - \frac{k}{m} }{ \frac{k}{m} } \right)
+ i \ln \Big(\sigma(1-q)^\ell\Big) 
+ \sum_{k = 1}^i 
\ln \,\psi\Big(\frac{k}{m}\Big)
 -\sum_{k = 1}^i 
\ln \, \phi\Big(\frac{k}{m}\Big)\, .
\end{align*}
The apparent singularity of~$\pi_m$ will be simplified later with~$\delta_m$. 
The first term can be written in a simpler form: 
\[
\sum_{k=1}^i \ln \bigg( \frac{ 1 - \frac{k}{m} }{ \frac{k}{m} } \bigg) 
= \ln \bigg(\frac{1}{i!}\prod_{k=1}^i (m-k)\bigg)
= \ln \binom{ m }{ i }+ \ln \Big(1-\frac{i}{m}\Big)\,,
\]
we thus obtain
\begin{multline}
\label{lnpi}
\ln \pi_i = 
\ln \binom{m}{i}+ \ln \Big(1-\frac{i}{m}\Big)+
i \ln \Big(\sigma(1-q)^\ell\Big)\\
+ \sum_{k = 1}^i 
\ln \,\psi\Big(\frac{k}{m}\Big)
 -\sum_{k = 1}^i 
\ln \, \phi\Big(\frac{k}{m}\Big)\, .
\end{multline}
Our goal is to estimate these quantities in the asymptotic regime~\eqref{regime}.
We will use a tricky comparison between a series and an integral to estimate the binomial coefficient and classical Taylor formulas to develop the sums. 

Let us first suppose that~$\sigma(1-q)^\ell-1$ does not tend towards~$\sigma-2$, the other case is a bit more complicated and will be treated in section~\ref{sigmaDeux}.
Since~$L_q$, defined in expression~\eqref{Lq}, is the difference between~$\sigma-2$ and~$\sigma(1-q)^\ell-1$, the quantity~$L_q$ stays far from~$0$ in the asymptotic regime. Thus there exists a positive number~$\eta$ such that
\begin{equation}
\label{eta}
|L_q|>\eta\,.
\end{equation}
Since~$\ell q \rightarrow a\neq 0$, there exists also a number~$\lambda \in \, ]0,1[$ such that
\begin{equation}
\label{lambda}
-1+\lambda < L_q\, .
\end{equation}
We will have to study two different cases for~$L_q$ and we will often need to treat them separately:

$\bullet$ The case when~$L_q$ is negative (corresponding to the case~$\sigma < 2$), in which we have
\[
-1+\lambda < L_q < -\eta\,.
\]
$\bullet$ The case when~$L_q$ is positive (corresponding to the case~$\sigma > 2$), in which we have
\begin{equation}
\label{Lqsigma}
\eta < L_q < \sigma-1\, .
\end{equation}
The last upper bound does not require any hypothesis at all, it follows from the expression~\eqref{Lq} of~$L_q$ and the fact that~$(1-q)^\ell$ is positive.

\subsection{Estimation of the binomial coefficient}
We now focus on the estimation of the binomial coefficient~$\ln \binom{m}{i}$.
A tricky comparison between a series and an integral derived by Robbins \cite{robbinsRemarkStirlingFormula1955a} yields the following inequalities:
\[
\forall n\geq 1 \qquad
\frac{1}{12n+1} < \ln n! -n\ln n +n -\frac{1}{2}\ln (2 \pi n)< \frac{1}{12n}\,.
\]
Let~$S(i)$ stand for
\begin{equation}
\label{S(i)}
S(i) = \ln \binom{m}{i} + m\Big(1-\frac{i}{m}\Big)\ln\Big(1-\frac{i}{m}\Big)+m\frac{i}{m}\ln\frac{i}{m} +\frac{1}{2}\ln\bigg(m\frac{i}{m}\Big(1-\frac{i}{m}\Big)\bigg)\,,
\end{equation}
computations then lead to, for every~$i \in \{1,\cdots, m-1\}\, ,$
\begin{multline}
\frac{1}{12m+1}-\frac{1}{12i}-\frac{1}{12(m-i)}\\
\leq
S(i) + \frac{1}{2} \ln (2\pi)
\leq\\
\frac{1}{12m}-\frac{1}{12i+1}-\frac{1}{12(m-i)+1}\,.
\end{multline}
we have the uniform bound: 
\begin{equation}
\label{ResteStirling}
\forall i\in \{1, \cdots, m-1\}\qquad
|S(i)| \leq \frac{1}{6}+\frac{1}{2}\ln (2\pi) \leq 2\,.
\end{equation}
The case~$i=m$ is easy since we have~$\ln\binom{m}{i}=0$ in this case.
\subsection{Expansion of the Riemann sums}
Let us now consider both the Riemann sum~$\sum \ln \psi$ and~$\sum \ln \phi$, we first write them as integrals.

For any function~$f$ of class~$C^2$ on~$[0, 1]$ and for~$k\in \{1, \cdots, m\}$, the Taylor-Lagrange formula applied to~$f$ between the points~$s\in[\frac{k-1}{m}, \frac{k}{m}]$ and~$\frac{k}{m}$ gives 
\[
\exists \eta_s^k\in \Big]s, \frac{k}{m} \Big[\qquad
f(s) = f\Big(\frac{k}{m}\Big) + \Big(s-\frac{k}{m}\Big)f'\Big(\frac{k}{m}\Big) + \Big(s-\frac{k}{m}\Big)^2\frac{f''(\eta_s^k)}{2}\,.
\]
Integrating this equality between the points~$\frac{k-1}{m}$ and~$\frac{k}{m}$ gives
\begin{equation}
\label{riem}
\int_{\frac{k-1}{m}}^{\frac{k}{m}}f(s)\, \diff s = \frac{1}{m} f\Big(\frac{k}{m}\Big)+ \Big(-\frac{1}{2m^2}\Big)f'\Big(\frac{k}{m}\Big) +\int_{\frac{k-1}{m}}^{\frac{k}{m}}\Big(s-\frac{k}{m}\Big)^2\frac{f''(\eta_s^k)}{2}\, \diff s\,.
\end{equation}
The last term will be negligible, we call it~$R_2$:
\[
R_2(k) = \int_{\frac{k-1}{m}}^{\frac{k}{m}}\Big(s-\frac{k}{m}\Big)^2\frac{f''(\eta_{s}^k)}{2}\, \diff s\,.
\]
We sum the expression~\eqref{riem} for~$k$ varying from~$1$ to~$i$ and we get
\begin{equation}
\label{Riemannf}
\int_{0}^{\frac{i}{m}}f(s) \diff s = \frac{1}{m} \sum_{k=1}^{i} f\Big(\frac{k}{m}\Big)- \frac{1}{2 m^2} \sum_{k=1}^{i} f'\Big(\frac{k}{m}\Big) +\sum_{k=1}^{i}R_2(k)\,.
\end{equation}
Similarly, we apply the Taylor-Lagrange formula at order~$1$ to~$f'$, we integrate and we sum to obtain
\begin{equation}
\label{Riemannf'}
\int_{0}^{\frac{i}{m}}f'(s) \diff s = \frac{1}{m} \sum_{k=1}^{i} f'\Big(\frac{k}{m}\Big)+ \sum_{k=1}^{i} \int_{\frac{k-1}{m}}^{\frac{k}{m}}\Big(s-\frac{k}{m}\Big)f''(\zeta_{s}^k)\, \diff s\, ,
\end{equation}
for some~$\zeta_{s}^k$ between~$s$ and~$\frac{k}{m}$.
We set
\begin{equation*}
R_1(k) =\int_{\frac{k-1}{m}}^{\frac{k}{m}}\Big(s-\frac{k}{m}\Big)f''(\zeta_{s}^k)\, \diff s\, .
\end{equation*}
We combine the two formulas~\eqref{Riemannf} and~\eqref{Riemannf'} and we obtain
\[
\sum_{k=1}^{i} f\Big(\frac{k}{m}\Big) =m\int_{0}^{\frac{i}{m}}f(s) \diff s + \frac{1}{2}\int_{0}^{\frac{i}{m}}f'(s) \diff s- \frac{1}{2}\sum_{k=1}^{i} R_1(k) -m\sum_{k=1}^{i}R_2(k)\, .
\]
Since~$f$ is a primitive of~$f'$, we get the following Lemma.
\vspace{10pt}
\begin{lemma}
\label{riemann}
For any function~$f$ of class~$C^2$ on~$[0, 1]$ and for every~$i\in\{1,\cdots, m\}$, we have
\[
\sum_{k=1}^{i} f\Big(\frac{k}{m}\Big) =
m\int_{0}^{\frac{i}{m}}f(s) \diff s + \frac{1}{2}\bigg(f \Big( \frac{i}{m} \Big) - f ( 0 )\bigg)- \frac{1}{2}\sum_{k=1}^{i} R_1(k) -m\sum_{k=1}^{i}R_2(k)\,,
\]
with, for~$k \in \{1, \cdots, i\}$,
\[R_1(k) =\int_{\frac{k-1}{m}}^{\frac{k}{m}}\Big(s-\frac{k}{m}\Big)f''(\zeta_{s}^k)\, \diff s\,,\]
\[R_2(k) = \int_{\frac{k-1}{m}}^{\frac{k}{m}}\Big(s-\frac{k}{m}\Big)^2\frac{f''(\eta_{s}^k)}{2}\, \diff s\,,\]
where~$\zeta_{s}^k$ and~$\eta_{s}^k$ belong to the interval~$]s,\frac{k}{m}[$.
\end{lemma}
We apply Lemma~\ref{riemann} to the function~$\ln \phi$. Let us set
\begin{equation}
\label{fphi}
f(x) = \ln \phi(x)= \ln \bigg( 1 - (1-q)^\ell Q^{\theta} +\Big( L_q + (1-q)^\ell Q^{\theta} \Big)x \bigg)\,.
\end{equation} 
The second derivative of~$f$ is
\[
f''(x) = \frac{ -\Big(L_q + (1-q)^\ell Q^{\theta}\Big)^2 }{ \Big(1-(1-q)^\ell Q^{\theta}+ (L_q+(1-q)^\ell Q^{\theta})x\Big)^2 }\, .
\]
We set \[a=\frac{-1+(1-q)^\ell Q^{\theta}}{L_q+(1-q)^\ell Q^{\theta}}\,,\] so we can write~$|f''|$ as

\[
|f''|(x) = \frac{ 1 }{ (x-a)^2 }\, .
\]
Note that~$q$ is small, as well as~$(1-q)^\ell Q^\theta$, and therefore~$a$ is not zero. 
Let us distinguish two cases:

$\bullet$ If~$L_q$ is negative, then~$a$ is positive, we even have~$a$ strictly greater than~$1$ because~$L_q > -1$. The second derivative of~$f$ is then uniformly bounded on the interval~$[0,1]$ by~$|f''|(1)$ which is equal to
\[
|f''|(1) = \bigg(\frac{ L_q + (1-q)^\ell Q^{\theta} }{ L_q +1}\bigg)^2\,.
\]
Since we have~$-1+\lambda < L_q < \sigma -1$ according to~\eqref{lambda}, we have that 
\[
|f''|(1) \leq \bigg(\frac{ L_q + (1-q)^\ell Q^{\theta} }{\lambda}\bigg)^2\leq \frac{\sigma^2}{\lambda^2}\,,
\]
and therefore the function~$f''$ is uniformly bounded on~$[0,1]$.

$\bullet$ If~$L_q$ is positive, then~$a$ is negative, and~$|f''|$ is decreasing, continuous on~$[0,1]$, thus~$|f''|$ is bounded by 
\[
|f''|(0) = \bigg(\frac{ L_q + (1-q)^\ell Q^{\theta} }{ 1-(1-q)^\ell Q^{\theta}}\bigg)^2 \leq (\sigma-1)^2 \leq 2 \sigma^2\,.
\]
Writing~$M$ for~$\max \Big(\, 2 \sigma^2, \,\frac{\sigma^2}{\lambda^2}\Big)$, we can bound both remainders as follows:
\[
\Big|\frac{1}{2}\sum_{k=1}^{i} R_1(k)\Big| \leq i\frac{M}{4 m^2} \leq \frac{M}{4 m} \,,
\] 
\[
\Big|m\sum_{k=1}^{i} R_2(k)\Big| \leq m i \frac{M}{6 m^3} \leq \frac{M}{6 m}\,.
\] 
In the end, Lemma~\ref{riemann} applied to the function~$\ln \phi$ gives
\begin{equation}
\label{sumphi}
\sum_{k=1}^{i} \ln \phi\Big(\frac{k}{m} \Big) =
m\int_{0}^{\frac{i}{m}} \ln \phi(s) \,\diff s 
+R^{(1)}\,,
\end{equation}
where
\[
R^{(1)} = \frac{1}{2}\bigg( \ln \phi\Big(\frac{i}{m} \Big) - \ln \phi(0)\bigg)-\frac{1}{2}\sum_{k=1}^{i} R_1(k) - m\sum_{k=1}^{i} R_2(k)\, .
\]
The quantity~$R^{(1)}$ is uniformly bounded by constants when~$m$ tends to infinity because~$\phi$ is a bounded function and the other terms tend to~$0$. 
We will treat in a similar way the function~$\ln \psi$, which we also call~$f$, 
\begin{equation}\label{fpsi}
f(x) = \ln \psi(x) = \ln \bigg( \frac{ Q^{\theta} }{ \sigma } +\Big( 1 - \frac{ Q^{\theta} }{ \sigma } \Big) x \bigg)\,.
\end{equation}
However, because of the very small value of~$\psi(0)$, the remainder term is more difficult to control here and we will have to be more precise in the upper bound of the second derivative of the function~$f$. The second derivative of the function~\eqref{fpsi} is
\[
f''(x)= \frac{ -\big(\sigma - Q^{\theta}\big)^2 }{ \big( Q^{\theta}+ (\sigma - Q^{\theta})x\big)^2 }\, .
\]
We define the associated remainders for~$k \in \{1, \cdots, i\}$:
\[R_1'(k) =\int_{\frac{k-1}{m}}^{\frac{k}{m}}\Big(s-\frac{k}{m}\Big)f''(\zeta_{s}^k)\, \diff s\,,\]
and
\[R_2'(k) = \int_{\frac{k-1}{m}}^{\frac{k}{m}}\Big(s-\frac{k}{m}\Big)^2\frac{f''(\eta_{s}^k)}{2}\, \diff s\,,\]
where~$\zeta_{s}^k$ and~$\eta_{s}^k$ belong to the interval~$]s,\frac{k}{m}[$.
We sum only for~$k$ between~$2$ and~$i$ in the formula~\eqref{riem} and we add~$f(1/m)$ on each side to get
\begin{equation}
\label{Riemannm}
\sum_{k=1}^{i} f\Big(\frac{k}{m}\Big) =m\int_{\frac{1}{m}}^{\frac{i}{m}}f(s) \diff s + \frac{1}{2}\bigg(f \Big( \frac{i}{m} \Big) + f \Big( \frac{1}{m} \Big)\bigg)- \frac{1}{2}\sum_{k=2}^{i} \Big(R_1'(k) -mR_2'(k)\Big)\,.
\end{equation}
We will bound the function~$f''$ on each sub-interval~$\big[\frac{k-1}{m}, \frac{k}{m}\big]$. Since the function~$|f''|$ is decreasing, then
\[
\sup_{[\frac{k-1}{m}, \frac{k}{m}]} |f''| \leq \frac{ \big(\sigma - Q^{\theta}\big)^2 }{ \big( Q^{\theta}+ (\sigma - Q^{\theta})\frac{k-1}{m}\big)^2 }\, .
\]
Therefore, 
\begin{align*}
\bigg|\sum_{k = 2}^{i} \int_{\frac{k-1}{m}}^{\frac{k}{m}}\Big(s-\frac{k}{m}\Big)f''(\zeta_{s})\, \diff s\bigg| 
&\,\leq\, \sum_{k = 2}^{i} \sup_{[\frac{k-1}{m}, \frac{k}{m}]} |f''| \int_{\frac{k-1}{m}}^{\frac{k}{m}}\Big|s-\frac{k}{m}\Big|\, \diff s \nonumber\\
&\leq \, \sum_{k = 2}^{i} \frac{ \big(\sigma - Q^{\theta}\big)^2 }{ \Big( Q^{\theta}+ (\sigma - Q^{\theta})\frac{k-1}{m}\Big)^2 } \frac{1}{2m^2}\, .
\end{align*}
Taking the first term out of the sum and shifting the indices, we obtain
\begin{align}
\label{R1f''}
\bigg|\sum_{k = 2}^{i} \int_{\frac{k-1}{m}}^{\frac{k}{m}}\Big(s-\frac{k}{m}\Big)f''(\zeta_{s})\, \diff s\bigg|
&\leq \frac{ \big(\sigma - Q^{\theta}\big)^2 }{ 2\Big( mQ^{\theta}+ \sigma - Q^{\theta}\Big)^2 } \\
& \qquad\qquad+ \frac{1}{2m^2}\sum_{k = 2}^{i-1} \frac{ \big(\sigma - Q^{\theta}\big)^2 }{ \Big( Q^{\theta}+ (\sigma - Q^{\theta})\frac{k}{m}\Big)^2 }\,.
\end{align}
We compare this second sum with an integral:
\[
\sum_{k = 2}^{i-1} \frac{ \big(\sigma - Q^{\theta}\big)^2 }{ \Big( Q^{\theta}+ (\sigma - Q^{\theta})\frac{k}{m}\Big)^2 } \leq 
m \int_{\frac{1}{m}}^{\frac{i-1}{m}}\frac{ \big(\sigma - Q^{\theta}\big)^2 }{ \Big( Q^{\theta}+ (\sigma - Q^{\theta})s\Big)^2 }\, \diff s\, .
\]
A change of variables gives
\[
\sum_{k = 2}^{i-1} \frac{ \big(\sigma - Q^{\theta}\big)^2 }{ \Big( Q^{\theta}+ (\sigma - Q^{\theta})\frac{k}{m}\Big)^2 }
\leq 
m \int_{Q^{\theta}+ (\sigma - Q^{\theta})\frac{1}{m}}^{Q^{\theta}+ (\sigma - Q^{\theta})\frac{i-1}{m}}\frac{ \sigma - Q^{\theta} }{ u^2 }\, \diff u\, .
\]
We calculate the integral and we obtain
\begin{align*}
\sum_{k = 2}^{i-1} \frac{ \big(\sigma - Q^{\theta}\big)^2 }{ \Big( Q^{\theta}+ (\sigma - Q^{\theta})\frac{k}{m}\Big)^2 }&\leq 
\frac{m \big(\sigma - Q^{\theta} \big)}{Q^{\theta}+ (\sigma - Q^{\theta})\frac{1}{m}} - \frac{m \big(\sigma - Q^{\theta} \big)}{Q^{\theta}+ (\sigma - Q^{\theta})\frac{i-1}{m}}\\
&\leq 
 \frac{m^2\big(\sigma - Q^{\theta}\big)}{mQ^{\theta}+ (\sigma - Q^{\theta})}\, ,
\end{align*}
where we removed the last negative term.
We finally bound~$\sum_{k = 2}^{i} R_1'(k)$, inequality~\eqref{R1f''} yields
\[
\bigg|\sum_{k = 2}^{i} \int_{\frac{k-1}{m}}^{\frac{k}{m}}\Big(s-\frac{k}{m}\Big)f''(\zeta_{s})\, \diff s\bigg| 
\leq 
\frac{ \big(\sigma - Q^{\theta}\big)^2 }{ 2\big( mQ^{\theta}+ \sigma - Q^{\theta}\big)^2 } + \frac{1}{2}\frac{\sigma - Q^{\theta}}{mQ^{\theta}+ (\sigma - Q^{\theta})}\, ,
\]
thus
\[
\Big|\sum_{k = 2}^{i} R_1'(k)\Big| \leq 1\,.
\]
The same work on~$R_2'(k)$ leads to
\[
\bigg|m\sum_{k = 2}^{i} \int_{\frac{k-1}{m}}^{\frac{k}{m}}\Big(s-\frac{k}{m}\Big)^2 \frac{f''(\zeta_{s})}{2}\, \diff s\bigg| 
\leq 
\frac{ \big(\sigma - Q^{\theta}\big)^2 }{ 6\big( mQ^{\theta}+ \sigma - Q^{\theta}\big)^2 } + \frac{1}{6}\frac{\sigma - Q^{\theta}}{mQ^{\theta}+ (\sigma - Q^{\theta})}\,,
\]
so we have
\[
\Big|m\sum_{k = 2}^{i} R_2'(k)\Big|\leq \frac{1}{3}\, .
\]
We gather these two remainders in 
\[
R^{(2)} =- \frac{1}{2}\sum_{k=2}^{i} R_1'(k) - m\sum_{k=2}^{i} R_2'(k)\,,
\]
and we get, according to~\eqref{Riemannm},
\begin{equation}
\label{sumpsi}
\sum_{k=1}^{i} \ln \psi\Big(\frac{k}{m}\Big) =
m\int_{1/m}^{i/m} \ln \psi(s) \diff s
+ \frac{1}{2}\bigg( \ln \psi\Big(\frac{i}{m}\Big) + \ln \psi\Big(\frac{1}{m}\Big)\bigg)+R^{(2)}\,,
\end{equation}
where~$|R^{(2)}|$ is uniformly bounded by~$1$.
\subsection{Computation of the sums}
We now compute the integrals,
the functions~$\psi$ and~$\phi$ are in fact affine functions, and the formula
\[
\int_{ 0 }^{ x } 
\ln(a+bs) \, \diff s =\frac{1}{ b} (a+bx)\ln (a+bx)-x - \frac{a}{ b}\ln a\,,
\]
will allow us to compute the two integrals appearing in equations~\eqref{sumphi} and~\eqref{sumpsi}. We have, according to the expressions~\eqref{psi} and~\eqref{phi} of the functions~$\psi$ and~$\phi$, 
\[
\int_{0}^{x} \ln \phi(s) \diff s = \frac{\phi(x) \ln \phi(x)}{L_q + (1-q)^\ell Q^{\theta}} - x -
\frac{\phi(0) \ln \phi(0)}{L_q + (1-q)^\ell Q^{\theta}}\,,
\]
\[
\int_{1/m}^{x} \ln \psi(s) \diff s = \frac{\sigma \psi(x) \ln \psi(x)}{\sigma -Q^{\theta} }- x -
\frac{\sigma \psi\Big(\frac{1}{m}\Big)\ln \psi\Big(\frac{1}{m}\Big) }{\sigma - Q^{\theta} } +\frac{1}{m}\,.
\]
Replacing the integral in formula~\eqref{sumphi} by the previous one gives
\begin{equation}
\label{sumlnphi}
\sum_{k = 1}^i \ln\phi\Big(\frac{k}{m}\Big) =
 \frac{m\phi\Big(\frac{i}{m}\Big) \ln \phi\Big(\frac{i}{m}\Big)}{L_q + (1-q)^\ell Q^{\theta}}- i 
-\frac{m\phi(0)\ln \phi(0) }{L_q + (1-q)^\ell Q^{\theta}} + R^{(1)}\,,
\end{equation}
where~$R^{(1)}$ is uniformly bounded by a constant. Formula~\eqref{sumpsi} also becomes
\begin{multline}
\label{sumlnpsi}
\sum_{k = 1}^i \ln \psi\Big(\frac{k}{m}\Big)= 
\frac{m\sigma \psi\Big(\frac{i}{m}\Big) \ln \psi\Big(\frac{i}{m}\Big)}{\sigma -Q^{\theta} }- i
-\frac{m\sigma \psi\Big(\frac{1}{m}\Big)\ln \psi\Big(\frac{1}{m}\Big) }{\sigma - Q^{\theta} } +1\\
+\frac{1}{2} \ln\psi\Big(\frac{i}{m}\Big)
+ \frac{1}{2} \ln \psi\Big(\frac{1}{m}\Big) + R^{(2)}\,,
\end{multline}
where~$R^{(2)}$ is uniformly bounded by constant terms. We will write~$R$ for the quantity~$-R^{(1)}+1+R^{(2)}$, there exists a constant~$C$ such that
\begin{equation}
\label{CsteRiemann}
 |R| \leq C\,.
\end{equation}
\subsection{First estimate of the persistence time}
We can now deal with the expectation of the persistence time, Lemma~\ref{vieMort} gives 
\[
E\Big(\tau_0 \, | \, N_0=1\Big)
=\sum_{i=1}^m\frac{1}{\delta_i}\exp\big(\ln \pi_i\big)\,.
\]
Replacing the sums in formula~\eqref{lnpi} with expressions~\eqref{sumlnphi} and~\eqref{sumlnpsi}, we formally obtain
\begin{multline}
\label{espetemp}
E\Big(\tau_0 \, | \, N_0=1\Big)=
\sum_{i=1}^m \frac{1-\frac{i}{m}}{\delta_i} \exp\Bigg(
\ln \binom{m}{i} + i\ln \Big(\sigma(1-q)^\ell\Big)\\
+\bigg(\frac{m\sigma \psi\Big(\frac{i}{m}\Big) }{\sigma -Q^{\theta} }+\frac{1}{2}\bigg)\ln \psi\Big(\frac{i}{m}\Big)
+\bigg(-\frac{m\sigma \psi\Big(\frac{1}{m}\Big) }{\sigma - Q^{\theta} } + \frac{1}{2} \bigg)\ln \psi\Big(\frac{1}{m}\Big)\\
-\frac{m\phi\Big(\frac{i}{m}\Big) \ln \phi\Big(\frac{i}{m}\Big)}{L_q + (1-q)^\ell Q^{\theta}}+ 
\frac{m\phi(0)\ln \phi(0) }{L_q + (1-q)^\ell Q^{\theta}} +R
\Bigg)\,.
\end{multline}
We begin by gathering the terms that do not depend on~$i$, we set
\[
K =
\exp\Bigg(
\bigg(-\frac{m\sigma \psi\Big(\frac{1}{m}\Big) }{\sigma - Q^{\theta} } + \frac{1}{2}\bigg)\ln \psi\Big(\frac{1}{m}\Big)+ 
\frac{m\phi(0)\ln \phi(0) }{L_q + (1-q)^\ell Q^{\theta}} \Bigg)\,.
\]
We replace the functions~$\psi$ and~$\phi$ by their expressions~\eqref{psi} and~\eqref{phi}, and we obtain
\begin{multline}
\label{K}
K=\exp\Bigg(
\bigg(-\frac{mQ^\theta}{\sigma-Q^\theta}-\frac{1}{2}\bigg)\ln\bigg(\frac{Q^\theta}{\sigma}+\Big(1-\frac{Q^\theta}{\sigma}\Big)\frac{1}{m}\bigg)\\
+m\frac{1-(1-q)^\ell Q^\theta}{L_q + (1-q)^\ell Q^{\theta}}\ln\left(1-(1-q)^\ell Q^\theta\right)\Bigg)\,.
\end{multline}
When we substitute~$\delta_i$ by its expression~\eqref{delta}, the expectation~\eqref{espetemp} becomes,
\begin{multline}
\label{sumAvantStirling}
E\Big(\tau_0 \, | \, N_0=1\Big)=K
\sum_{i=1}^m \frac{\sigma\frac{i}{m}+1-\frac{i}{m}}{\sigma(1-q)^\ell}
 \exp\Bigg(
\ln \binom{m}{i} + i\ln \Big(\sigma(1-q)^\ell\Big)\\
+\bigg(\frac{m\sigma \psi\Big(\frac{i}{m}\Big)}{\sigma -Q^{\theta} } -\frac{1}{2}\bigg) \ln \psi\Big(\frac{i}{m}\Big)
-\frac{m\phi\Big(\frac{i}{m}\Big) \ln \phi\Big(\frac{i}{m}\Big)}{L_q + (1-q)^\ell Q^{\theta}} +R
\Bigg)\,.
\end{multline}
In the sequel, we will sum only to~$m-1$ and let the last term out.
In the main sum, we introduce the expansion~\eqref{S(i)} of the binomial coefficient.
We still denote by~$R$ the remainder term when we take into account the new quantity~$S$ appearing in~\eqref{S(i)}. Thanks to the bound~\eqref{ResteStirling}, we still have
\begin{equation}
\label{RC}
 |R|\leq C'\,,
\end{equation}
 for some constant~$C'$. The last logarithmic term of~\eqref{ResteStirling} is not defined when~$i=m$, so from now onwards we will isolate this term from the sum.

We introduce some new notations. First let us set the function~$\widetilde{F}$ to write in a simpler way the term that contain~$\phi$:
\begin{equation}
\label{FtildeDef}
 \widetilde{ F}(x) = \frac{ 1 + L_q x }{ L_q }\ln (1 + L_q x ) 
 -\frac{\phi(x) \ln \phi(x)}{L_q + (1-q)^\ell Q^{\theta}}\,.
\end{equation}
This function allows us to write
\begin{equation*}
-\frac{\phi(x) \ln \phi(x)}{L_q + (1-q)^\ell Q^{\theta}}
=
- \frac{ 1 + L_q x }{ L_q } \ln (1 + L_q x ) + \widetilde{ F}(x)\,.
\end{equation*}
We then define the function~$F$ by gathering the terms that are multiplied by~$m$ in the exponential:
\begin{equation}
\label{F}
F ( x ) = 
- \left( 1 - x \right) \ln \left( 1 - x \right) + x \ln \Big(\sigma(1-q)^\ell\Big) 
- \frac{ 1 + L_q x }{ L_q } \ln (1 + L_q x )\,.
\end{equation}
To express the sum~\eqref{sumAvantStirling}, we first replace~$\psi$ by its expression~\eqref{psi} and we write
\[
\frac{m\sigma \psi\big(\frac{i}{m}\big)}{\sigma -Q^{\theta} } = \frac{mQ^\theta}{\sigma -Q^{\theta} } + m\frac{i}{m}\,.
\]
We gather together the two quantities~$ m\frac{i}{m}\ln \psi\big(\frac{i}{m}\big)$, and~$-m\frac{i}{m}\ln\frac{i}{m}$ from the binomial expansion, we write them as 
\[
m\frac{i}{m}\ln \psi\Big(\frac{i}{m}\Big)-m\frac{i}{m}\ln\frac{i}{m} = m \frac{i}{m}\ln\bigg(\frac{Q^\theta}{\sigma \frac{i}{m}}+1-\frac{Q^\theta}{\sigma}\bigg)\,.
\]
Finally, we put the remaining terms together in a function~$G$:
\begin{multline}
\label{G}
G ( x ) =
\ln \frac{ (\sigma-1) x + 1 }{ \sigma (1-q)^{\ell}}
+ \Big(\frac{mQ^\theta}{\sigma-Q^\theta}-\frac{1}{2}\Big)\ln\bigg(\frac{Q^\theta}{\sigma}+\Big(1-\frac{Q^\theta}{\sigma}\Big)x\bigg)\\
+ m x\ln\bigg(\frac{Q^\theta}{\sigma x}+1-\frac{Q^\theta}{\sigma}\bigg)
-\frac{1}{2}\ln\Big( mx(1-x)\Big) +m\widetilde{F}(x)\,.
\end{multline}

With these notations, formula~\eqref{sumAvantStirling} reduces to
\begin{equation}
\label{FFG}
E\Big(\tau_0 \, | \, N_0=1\Big)=K\sum^{m-1}_{i=1} \, \exp \bigg( m F \Big(\frac{i}{m}\Big) +G \Big( \frac{i}{m} \Big) +R
\bigg)
+T,
\end{equation}
where~$T$ is the term of index~$m$:
\begin{align*}
T&= \frac{K}{(1-q)^\ell }\exp\bigg(m\ln\Big(\sigma (1-q)^\ell\Big)-m\frac{\big(1+L_q\big)\ln\big(1+L_q\big)}{L_q +(1-q)^\ell Q^\theta}+R\bigg)\\
&=K\exp\Big(mF(1)+ m\widetilde{F}(1)+C''\Big)\,,
\end{align*}
where we used the facts that~$\phi(1) = 1+L_q$, that~$\psi(1) = 1$, and that~$(1-q)^\ell$ does not go to zero because~$\ell q$ does not go to infinity.

The function~$F$ in the sum~\eqref{FFG} is multiplied by~$m$, which tends towards infinity. If we can bound the function~$G$ uniformly over the interval~$[\frac{1}{m}, 1-\frac{1}{m}]$, the indices around the maximum of the function~$F$ will govern the asymptotic behavior of the sum. 
\subsection{Function~$G$}
Let us study the terms that appear in the expression~\eqref{G} of the function~$G$. Let~$x$ belong to the interval~$[\frac{1}{m}, 1-\frac{1}{m}]$.

$\bullet$ For the third term, we have, since~$\sigma > 1$, 
\begin{equation*}
\left|m x \ln\bigg( \frac{ Q^{\theta} }{ \sigma x} + 1 - \frac{ Q^{\theta} }{ \sigma } \bigg)\right| 
\leq mx \frac{Q^\theta}{\sigma}\Big(\frac{1}{x}-1\Big)
\leq m Q^\theta\,.
\end{equation*}
$\bullet$ The second term gives, since~$x > 1/m$,
\[
\left| 
\ln\bigg( \frac{ Q^{\theta} }{ \sigma } + \Big(1 - \frac{ Q^{\theta} }{ \sigma } \Big) x \bigg)\right| \leq 
-\ln\bigg( \frac{ Q^{\theta} }{ \sigma } + \Big(1 - \frac{ Q^{\theta} }{ \sigma } \Big) \frac{1}{m}\bigg)\,.
\]
Since the function~$-\ln$ is convex, then
\[
-\ln\bigg( \frac{ Q^{\theta} }{ \sigma } + \Big(1 - \frac{ Q^{\theta} }{ \sigma } \Big) \frac{1}{m}\bigg)\leq
\Big(1 - \frac{ Q^{\theta} }{ \sigma } \Big)\Big(-\ln \frac{1}{m}\Big) \leq \ln m\,.
\]
Since~$\sigma>1$ and~$Q$ tends towards~$0$, we get, for every~$x$ in~$[\frac{1}{m}, 1-\frac{1}{m}]$,
\begin{equation}
\label{M(theta)}
 \bigg|\Big(\frac{mQ^\theta}{\sigma-Q^\theta}-\frac{1}{2}\Big)\ln\bigg(\frac{Q^\theta}{\sigma}+\Big(1-\frac{Q^\theta}{\sigma}\Big)x\bigg)\bigg|
\leq (mQ^\theta+1)\ln m\,.
\end{equation}
$\bullet$ Since the function~$x\mapsto x(1-x)$ is always smaller than~$1/4$, then the last term of function~$G$ can be controlled by
\begin{equation*}
\left| 
-\frac{1}{2}\ln\Big( mx(1-x)\Big)\right| 
\leq \frac{1}{2}\ln\Big( \frac{m}{4}\Big)
\leq \ln m\,.
\end{equation*}
$\bullet$ We now bound from above the expression~\eqref{FtildeDef} of~$\widetilde{F}$. According to the expression~\eqref{phi} of~$\phi$, 
\[
\phi(x) = 1+L_q x -(1-q)^\ell Q^\theta(1-x)\,,
\]
so
\[
\ln \phi(x) =\ln \big(1 + L_q x\big) + \ln \bigg(1 -\frac{(1-q)^\ell Q^{\theta} (1-x)}{1+L_qx} \bigg).
\]
A triangular inequality gives 
\begin{multline}
\label{MajTempFtil}
|\widetilde{F}(x)| \leq 
\Big|\big(1 + L_q x\big)\ln (1+L_q x )\Big|\bigg| \frac{1}{L_q}- \frac{1}{L_q+(1-q)^\ell Q^{\theta} }\bigg| \\
+ \Big|\frac{ 1 + L_q x }{ L_q +(1-q)^\ell Q^{\theta}} \Big|\Big|\ln \bigg(1 -\frac{(1-q)^\ell Q^{\theta} (1-x)}{1+L_qx} \bigg)\Big|\\
+\Big|\frac{(1-q)^\ell Q^{\theta}(1-x) }{ L_q + (1-q)^\ell Q^{\theta} }\Big| \Big|\ln\phi(x)\Big|\,.
 \end{multline}
$\bullet$ We compute the difference of the two fractions in the first term and we get, when~$Q$ is small enough, and with the help of the inequality~\eqref{eta}
\[
\bigg| \frac{1}{L_q}- \frac{1}{L_q+(1-q)^\ell Q^{\theta} }\bigg| \leq \frac{(1-q)^\ell Q^{\theta}}{|L_q||L_q+(1-q)^\ell Q^{\theta}|}\leq 
\frac{ Q^{\theta}}{\eta^2}\,.
\]
$\bullet$ For the second term, since~$L_q > -1$, we have~$ 1+L_q x > 1-x$, so that
\begin{equation}
\label{lncomplique}
0\leq -\ln \bigg(1 -\frac{(1-q)^\ell Q^{\theta}(1-x) }{ 1+L_q x} \bigg)\leq -\ln \big(1 -Q^{\theta} \big)\,.
\end{equation}
Expression~\eqref{Lq} gives a simple bound on~$|1+L_q x|$:
\begin{equation}
\label{1+Lqx}
|1+L_q x| \leq 1+\sigma x\leq 2 \sigma\,.
\end{equation}

Since~$\phi$ is smaller than~$1+L_q x$, the bound~\eqref{MajTempFtil} gives thanks to the three previous inequalities, for~$x$ in~$[0, 1]$,
\[
|\widetilde{F}(x)| \leq 
\frac{ Q^{\theta}\, 2 \sigma|\ln (1+L_q x )|}{\eta^2} 
+ \frac{ 2 \sigma }{ \eta} \Big|\ln \big(1 -Q^{\theta} \big)\Big|
+\frac{Q^{\theta}}{ \eta } \Big|\ln(1+L_q x)\Big|\,.
\]
The inequality~$-\ln(1-u)\leq \ln(1+2u)$ holds as soon as~$u\leq 1/2$. For~$Q$ small enough, we have therefore
\[
\Big|\ln \big(1 -Q^{\theta} \big)\Big| \leq 
2Q^\theta\,.
\]
We now deal with~$\ln(1+L_q x)$:

$\bullet$ if~$L_q$ is positive, we can use bound~\eqref{Lqsigma} to bound from above by
\[
\ln(1+L_q x)\leq \ln \sigma\,.
\]
$\bullet$ If~$L_q$ is negative, the function~$x\mapsto 1+L_q x$ decreases with~$x$, so according to~\eqref{lambda},
\[
1+L_q x \geq 1+L_q \geq \lambda\,.
\]
In both cases, we have 
\[
|\ln(1+L_q x)| \leq \ln(1/\lambda)\,.
\]
We choose~$1/\lambda$ to be greater than~$\sigma$, and we get
\[
|\widetilde{F}(x)| \leq 
\frac{ Q^{\theta}\, 2 \sigma\ln (1/\lambda)}{\eta^2} 
+ \frac{ 4 \sigma Q^\theta }{ \eta}
+\frac{Q^{\theta}\ln(1/\lambda)}{ \eta } \,.
\]
Therefore
\begin{equation}
\label{FtildeMaj}
\sup_{[0,1]} |\widetilde{F}| \leq \frac{8\sigma \ln (1/\lambda)}{\eta^2} \,Q^\theta\,.
\end{equation}
Putting together formulas~\eqref{G},~\eqref{M(theta)}, and~\eqref{FtildeMaj}, we get for every~$x$ in the interval~$[\frac{1}{m},1-\frac{1}{m}]$,
\[
\big|G ( x )\big| \leq
\ln \frac{ 1}{(1-q)^{\ell}} + (mQ^\theta+1) \ln m\\
+ m Q^\theta
+\ln m + \frac{8\sigma \ln (1/\lambda)}{\eta^2} \,mQ^\theta\,.
\]
In the asymptotic regime~\eqref{regime}, we have thus
\begin{equation}
\label{MajG}
\sup_{[\frac{1}{m}, \frac{m-1}{m}]} \left|G(x)\right| \, \leq \, 2(1+mQ^\theta)\ln m\,.
\end{equation}
We will use this upper bound several times in the sequel.
\subsection{Function~$F$}
We are now looking for the maximum of~$F$ on~$[0,1]$, which we call~$\rho^*$. The function~$F$ defined in expression~\eqref{F} admits for first derivative
\[
F'(x)
=
\ln (1-x) + 
\ln \Big(\sigma(1-q)^\ell\Big) - 
\ln \big(1+ L_q x\big)\,,
\]
and for second derivative
\begin{equation}
\label{F''}
F''(x)=-\frac{1}{1-x} - \frac{L_q}{1+ L_q x}\,.
\end{equation}
According to the assumption~\eqref{lambda}, we have
\begin{equation}
\label{F''neg}
F''(x)\leq -\frac{1}{1-x} - \frac{L_q}{1+ L_q x} \leq -\frac{1+L_q}{1+ L_q x} \leq -\lambda <0\,.
\end{equation}
The function~$F$ is therefore concave and its unique critical point~$\rho^*$ is a maximum.
This point satisfies the following equation:
\begin{equation}
\label{equationRho}
\ln(1-\rho^*)+
\ln \Big(\sigma(1-q)^\ell\Big) -
\ln \big(1+ 
L_q \rho^*\big) 
= 0\,.
\end{equation}
Therefore, we have
\[
\rho^* 
= 
\frac{ \sigma (1-q)^\ell -1 }{ \sigma - 1 }\,.
\]
Under the condition~\eqref{regimeconditionrho}, the product~$\sigma (1-q)^\ell$ is asymptotically greater than~$1$, so~$\rho^*$ is indeed positive.
Moreover, according to expression~\eqref{Lq} of~$L_q$, we can write
\[
\rho^*=\frac{\sigma-2-L_q}{\sigma-1}\,.
\]
The quantity~$1+L_q \rho^*$ will appear often in the sequel, the equation~\eqref{equationRho} provides the following expression for this quantity:
\begin{equation}
\label{unplusLqrho}
1+L_q \rho^*
=\, \sigma(1-q)^\ell (1-\rho^*)\,,
\end{equation}
and since~$L_q>-1+\lambda$, 
\begin{equation}
\label{rhopetit}
\rho^* < 1 - \frac{\lambda}{\sigma-1}\,.
\end{equation}
We will replace the term~$F(i/m)$ in the sum~\eqref{FFG} by its Taylor development around~$\rho^*$, so we first calculate~$F(\rho^*)$:
\[
F(\rho^*)
=\,
-(1-\rho^*)\ln(1-\rho^*)+\rho^*\ln\Big(\sigma(1-q)^\ell\Big) 
-\Big(\frac{1}{L_q }+\rho^*\Big)\ln \left(1+L_q\rho^*\right)\,.
\]
Equation~\eqref{unplusLqrho} gives
\[
\Big(\frac{1}{L_q }+\rho^*\Big)\ln \left(1+L_q\rho^*\right)=\Big(\frac{1}{L_q }+\rho^*\Big)\ln \left(\sigma(1-q)^\ell(1-\rho^*)\right)\, .
\]
Splitting the logarithmic term, we obtain
\[
F(\rho^*)
=\, -\frac{1+L_q}{L_q}\ln (1-\rho^*)-\frac{1}{L_q }\ln \left(\sigma(1-q)^\ell\right) .
\]
Replacing~$L_q$ with its expression~\eqref{Lq} finally gives
\begin{equation}
\label{F(rho)}
F(\rho^*) = \,
\frac
{
\sigma(1-(1-q)^\ell)
\ln \frac{ \sigma(1-(1-q)^\ell) }{ \sigma-1} + \ln \Big(\sigma (1-q)^\ell\Big)
}{
1-\sigma\big(1-(1-q)^\ell\big)
}\,.
\end{equation}
This quantity will often appear in the sequel, we define the function~$\varphi$ as
\begin{equation}
\label{varphiDef}
\varphi(x) =\frac{
\sigma(1-x)
\ln \frac{ \sigma(1-x) }{ \sigma-1} + \ln (\sigma x)
}{
1-\sigma(1-x)
}\,,
\end{equation}
so that we can write
\[
F(\rho^*) = \varphi\Big((1-q)^\ell\Big)\,.
\]
We will also need~$F''(\rho^*)$, according to~\eqref{unplusLqrho}, we have
\begin{align*}
F''(\rho^*) &=
-\frac{1 }{ 1-\rho^* }
-\frac{L_q }
{1+L_q \rho^*}\\
&=
\Big(-1 
-\frac{L_q }
{\sigma(1-q)^\ell }\Big)\frac{1}{1-\rho^*}\, ,
\end{align*}
so 
\begin{equation}
\label{F''(rho)}
F''(\rho^*)=
- \frac{
\left(\sigma-1\right)^2
}{
\sigma^2(1-q)^\ell (1-(1-q)^\ell)
}\,.
\end{equation}
\section{Implementation of Laplace's method}
We now introduce a notation for the sum~\eqref{FFG}: we set
\begin{equation}
\label{Sm}
S_m = \sum^{m-1}_{i=1} \, \exp \bigg( m F \Big(\frac{i}{m}\Big) + G\Big( \frac{i}{m} \Big)\bigg)\,.
\end{equation}
The expectation of the persistence time can thus be written as
\begin{equation}
\label{ES}
E\Big(\tau_0 \, | \, N_0=1\Big) =KS_m\,e^{R_1}+R_2\,,
\end{equation}
where, according to expression~\eqref{RC},
\[
|R_1| \leq |C''| \,,
\]
and
\[
 R_2 = K \exp\Big(mF(1)+ m\widetilde{F}(1)+C''\Big)\,.
\]
Our next objective is now to estimate the sum~\eqref{Sm}. 
The main contributions in the sum will arise from terms whose indices lie around~$m\rho^*$, therefore, we will first estimate the sum truncated on a certain neighborhood of~$m\rho^*$.We take
\begin{equation}
\label{deltam}
\delta = m^{2/3},
\end{equation}
and we set~$[i_-, i_+]$ to be the interval on which we will sum, where
\begin{equation*}
i_- = \max\Big(\lfloor m\rho^* -\delta\rfloor, \,0 \Big)+1\,,
\end{equation*}
\begin{equation*}
i_+ = \lfloor m\rho^* +\delta\rfloor\,.
\end{equation*}
Since~$i_- \geq 1$ and~$\rho^*<1-\frac{\lambda}{\sigma-1}$ according to inequality~\eqref{rhopetit}, the interval~$[i_-, i_+]$ is strictly included in~$[1,m-1]$. 
Our goal is now to estimate the sum
\begin{equation}
\label{Sm(delta)}
S_m(\delta) = \sum_{i=\, i_-}^{ i_+} \, \exp \bigg( m F \Big(\frac{i}{m}\Big) + G \Big( \frac{i}{m} \Big) \bigg)\,.
\end{equation}
Recalling that~$\rho^*$ is the maximum of the function~$F$, this quantity is related to the expression~\eqref{Sm} through the inequalities
\[
 S_m(\delta)\leq
S_m
\leq S_m(\delta)
+m \exp\Big( m F (\rho^*)+\sup_{[\frac{1}{m}, \frac{m-1}{m}]} G\Big)\,.
\]
We thus obtain, according to the inequality~\eqref{MajG},
\begin{equation}
\label{SmSm(delta)}
S_m = S_m(\delta) + R_3\,,
\end{equation}
with
\[
0\leq R_3 \leq \exp\Big( m F (\rho^*)
+2(1+mQ^\theta)\ln m\Big)\,.
\]
The Taylor-Lagrange formula at order~$3$ for~$F$ allows us to estimate the expression~\eqref{Sm(delta)} of~$S_m(\delta)$:
\begin{multline}
\label{SmTaylor}
 S_m(\delta)
 =\sum_{i=\,i_-}^{i_+} \exp \Bigg( m \bigg(F (\rho^*)+\Big(\frac{i}{m}-\rho^*\Big)^2\frac{F'' \left(\rho^*\right)}{2}+\Big(\frac{i}{m}-\rho^*\Big)^3\frac{F''' \left(\eta_i\right)}{6}\bigg)\\
 + G \left(\frac{i}{m}\right) \Bigg)\,,
\end{multline}
where~$\eta_i$,~$i_- \leq i \leq i_+$ is a real number between~$i_-/m$ and~$i_+/m\,.$
\subsection{Control of the remainders}
By differentiating the expression~\eqref{F''} of~$F''$, we get successively
\[
F'''(x) = -\frac{1}{(1-x)^2}+\frac{L_q^2 }{(1+L_q x)^2}\,,
\]
and
\[
F^{(4)}(x) =-\frac{2}{(1-x)^3}-\frac{2 L_q^3 }{(1+L_q x)^3}\,.
\]
As~$L_q > -1$, we have~$-2L_q^3< 2$, and~$0< 1-x \leq 1+L_qx$, which leads to
\[
\forall x \in [0,1]\qquad F^{(4)}(x) < -2\bigg(\frac{1}{(1-x)^3}-\frac{1}{(1+L_q x)^3}\bigg)< 0\,.
\]
The function~$F^{(4)}$ is negative on~$[0,1]$. On any compact interval included in~$[0,1]$, the function~$F'''$ is uniformly bounded by its values at the boundaries. Since~$i_- > 0$ and~$i_+ \leq m\rho^* +\delta$, then
\[
\forall x \in \Big[0, \rho^* +\frac{\delta}{m}\Big]\qquad F'''\Big(\rho^* +\frac{\delta}{m}\Big) \leq F'''(x) \leq F'''(0)\,.
\]
For the lower bound, we bound the positive term by~$0$ and we obtain
\[
F'''\Big(\rho^* +\frac{\delta}{m}\Big)\geq 
-\frac{1}{(1-\rho^* -\frac{\delta}{m})^2}+\frac{L_q^2 }{(1+L_q \rho^* +L_q\frac{\delta}{m})^2}\geq -\Big(\frac{\sigma-1}{\lambda}\Big)^2\,,
\]
according to inequality~\eqref{rhopetit}.
Therefore, on the interval~$[i_-/m, i_+/m]$, we have 
\[
-\Big(\frac{\sigma}{\lambda}\Big)^2\leq F''' \leq\sigma^2\,,
\]
thus, the expression~$ m\big(\frac{i}{m}-\rho^*\big)^3 F'''(\eta_i)$ is uniformly bounded by constants for every~$i$ in the interval~$[i_-, i_+]$, we write
\begin{equation}
\label{MajF'''}
\forall i \in [i_-, i_+]\qquad
\Big|m\Big(\frac{i}{m}-\rho^*\Big)^3 F'''(\eta_i)\Big| \leq \frac{\sigma^2}{\lambda^2}\,.
\end{equation}

\subsection{The main term}
Let~$T_m(\delta)$ be the sum we want to estimate, 
\begin{equation}
\label{Tm(delta)}
T_m(\delta) = \sum_{i=i_-}^{i_+} \exp \bigg( m\Big(\frac{i}{m}-\rho^*\Big)^2\frac{F'' \left(\rho^*\right)}{2}\bigg) \,.
\end{equation}
From the expression~\eqref{SmTaylor}, we have
\begin{equation}
\label{Smdvpe}
S_m(\delta)=\exp \Big( m F \left(\rho^*\right)\Big) \,T_m(\delta)\, e^{R_4},
\end{equation}
where, according to~\eqref{MajF'''} and~\eqref{MajG}, the remainder term~$R_4$ satisfies
\[
\big|R_4 \big| \leq \frac{\sigma^2}{\lambda^2}+2(1+mQ^\theta)\ln m \,.
\]
We will only need a rough approximation on~$T_m(\delta)$. Since~$F''(\rho^*)$ is negative, we first notice that 
\[
T_m(\delta) \leq m\,.
\]
We can also bound from below~$T_m(\delta)$ by one of the terms: for example, the term of index~$1+\lfloor m\rho^*\rfloor$
\[
T_m(\delta)\geq \exp\bigg(m\Big(\frac{1+\lfloor m\rho^*\rfloor}{m}-\rho^*\Big)^2\frac{F''(\rho^*)}{2}\bigg)\,.
\]
However, we have that 
\[
m\rho^*-1 \leq \lfloor m\rho^* \rfloor \leq m\rho^*\,,
\]
so
\[
T_m(\delta)\geq \exp\Big(\frac{F''(\rho^*)}{2m}\Big)\,
\]
and this bound goes to 1, since~$F''(\rho^*)$ is finite as we saw in expression~\eqref{F''(rho)} and in the hypothesis that the product~$\ell q$ does not tend towards~$0$ or~$\infty$.

%
Therefore, we have for~$m$ large enough
\[
\frac{1}{2}\leq T_m(\delta)\leq m\,.
\]
With the formula~\eqref{Smdvpe}, we rewrite the sum~\eqref{Sm(delta)} as
\begin{equation}
\label{Sm(delta)fin}
S_m(\delta) =\,\exp\Big( m F (\rho^*)\Big)\, R_5 \, e^{R_4}\,,
\end{equation}
where 
\[
\frac{1}{2} \leq R_5\leq m\,.
\]

We can now come back to the persistence time: formulas~\eqref{ES}, and~\eqref{SmSm(delta)} give
\[
E\Big(\tau_0 \, \big| \, N_0=1\Big) 
=K\Big(S_m(\delta)+R_3\Big)e^{R_1}+R_2\,.
\]
Replacing~$S_m(\delta)$ with its expression~\eqref{Sm(delta)fin}, we get 
\[
E\Big(\tau_0 \, | \, N_0=1\Big)
= K \exp\Big(m F (\rho^*)\Big)R_5 e^{R_4+R_1}
+ KR_3e^{R_1}
+R_2\, .
\]
We will show that the persistence time is close to the first term in this last expression, so we divide the equality by~$K \exp \big( m F \left(\rho^*\right)\big)$ and we look at the remainders.

$\bullet$ For the first one, we have, according to formulas~\eqref{SmSm(delta)} and~\eqref{ES},
\[
\bigg|\frac{KR_3e^{R_1}}{K \exp \big( m F \left(\rho^*\right)\big)}\bigg| \leq
\exp\Big(2(1+mQ^\theta)\ln m +C'' \Big)\,.
\] 

$\bullet$ For the second and last remainder term, we have, according to formula~\eqref{ES},
\[
\bigg|\frac{R_2}{ K\exp \big( m F \left(\rho^*\right)\big)} \bigg|\leq
\exp\Big(
-m\big(F(\rho^*)-F(1)\big)+m\widetilde{F}(1)+C'' \Big)\,.
\]
Since~$\rho^*$ is the maximum of function~$F$, and according to~\eqref{FtildeMaj}
\begin{align*}
\bigg|\frac{R_2}{ K\exp \big( m F \left(\rho^*\right)\big)} \bigg|&\leq
\exp\Big(
m\sup_{[0,1]} \widetilde{F}+C'' \Big) \\
&\leq 
\exp\Big(
m\frac{8\sigma \ln (1/\lambda)}{\eta^2} \,Q^\theta+C'' \Big)\,.
\end{align*}

%
%
%
%
We have then
\begin{equation}
\label{grosseExp}
E\Big(\tau_0 \, | \, N_0=1\Big) =
K\exp \bigg( m F \left(\rho^*\right)+O\Big( (1+mQ^\theta)\ln m \Big)\bigg)\,.
\end{equation}
We now develop this expression, remember that expression~\eqref{K} of~$K$ yielded
\begin{multline*}
K 
=\exp\Bigg(
\Big(\frac{mQ^\theta}{\sigma-Q^\theta}+\frac{1}{2}\Big)\bigg(\ln m -\ln\Big(\frac{mQ^\theta}{\sigma}+1-\frac{Q^\theta}{\sigma}\Big)\bigg) \\ 
+m\frac{1-(1-q)^\ell Q^\theta}{L_q + (1-q)^\ell Q^{\theta}}\ln\left(1-(1-q)^\ell Q^\theta\right)\Bigg)\,.
\end{multline*}
Some terms in~$K$ are of order smaller than~$(1+mQ^\theta)\ln m$, we include them in the remainder term, and we get according to definition~\eqref{varphiDef},
 \begin{multline*}
E\Big(\tau_0 \, | \, N_0=1\Big)=
 \exp \Bigg( 
m\varphi\Big((1-q)^\ell\Big) 
-\Big(\frac{mQ^\theta}{\sigma-Q^\theta}+\frac{1}{2}\Big)\ln\bigg(\frac{mQ^\theta}{\sigma}+1-\frac{Q^\theta}{\sigma}\bigg) \\
+ O\Big((1+mQ^\theta)\ln m\Big)\Bigg)\,.
\end{multline*}
Let us now focus on the complicated term, we have that
\[
\sigma-Q^\theta\geq \frac{\sigma}{2}\geq \frac{1}{2}\,,
\]
moreover, is~$m$ is large enough,
\[
1+\frac{mQ^\theta}{\sigma}-\frac{Q^\theta}{\sigma}\leq m\,.
\]
Therefore,
\[
0\leq \Big(\frac{mQ^\theta}{\sigma-Q^\theta}+\frac{1}{2}\Big)\ln\bigg(\frac{mQ^\theta}{\sigma}+1-\frac{Q^\theta}{\sigma}\bigg) 
\leq (1+mQ^\theta)\ln m\,
\]
so this term is at most of the order of the remainder term.
Finally, we obtain
\[
E\Big(\tau_0 \, | \, N_0=1\Big)=
 \exp \Bigg( 
m\varphi\Big((1-q)^\ell\Big)
+ O\Big((1+mQ^\theta)\ln m\Big)\Bigg)\,.
\]
Thus theorem~\ref{theoremTemps} is proved in the case where~$\sigma(1-q)^\ell-1$ does not tend to~$\sigma-2$.
\section{The error threshold}
Let us remind that the discovery time~$\tau^*$ is of order
\[
E\big(\tau^* \big)=
 \exp \Big( 
\ell \ln \kappa
+ o(\ell)\Big)\,.
\]
We are going to distinguish~$3$ cases according to the behaviour of the ratio~$m/\ell$.
\subsection{In the case~$m/\ell\rightarrow \infty$}
This case corresponds to a large population, however, we must ensure that the genome length~$\ell$ is not too small, because the comparison of the two times will only be relevant if the remainder terms are indeed smaller than the main terms.
Let us thus suppose that
\[
 \frac{\ell^2}{m\ln m}\rightarrow \infty\,.
\]
The condition we just set implies that 
\[
\frac{\ell}{(1+mQ^\theta)\ln m}\rightarrow \infty\,.
\]
%
%
We deduce that for the two times~$\tau_0$ and~$\tau^*$ to be of the same order, it is necessary that the ratio 
$$
\frac{m\varphi\Big((1-q)^\ell\Big)}{\ell \ln \kappa}
$$
stays bounded.
Therefore the quantity~$\varphi\big((1-q)^\ell\big)$ must go to~$0$.
Since the function~$\varphi$ is a bijection from~$[1/\sigma, 1]$ to~$[0, \ln\sigma]$, it implies that the quantity~$\sigma(1-q)^\ell-1$ must go to~$0$.
If it not the case, the time~$\tau^*$ will always be greater than the persistence time~$\tau_0$, thus the concentration of master sequences will be negligible. 

Let us now suppose that
\[
\sigma\neq2\,,
\]
and 
\[
\sigma(1-q)^\ell -1\rightarrow 0\,,
\]
and let us study the asymptotic behavior of the term~$ \varphi\big((1-q)^\ell\big)$.
Its asymptotic will be given by the first non-zero derivative of~$\varphi$ at point~$1/\sigma$. In order to expand~$\varphi$ around~$1/\sigma$, we write
\[
\varphi(x) = \varphi \bigg(\Big(x-\frac{1}{\sigma}\Big)+\frac{1}{\sigma}\bigg)\,,
\]
and we use the expression~\eqref{varphiDef} of~$\varphi$ to get
\[
\varphi(x) = \frac{\Big(\sigma-(\sigma x-1) -1\Big)
\ln \frac{ \sigma-(\sigma x-1) -1 }{ \sigma-1} + \ln \Big((\sigma x-1) +1\Big)
}{
1-\sigma(1-x)
}\,.
\]
We develop the expression in powers of~$\sigma x-1$ and we get
\begin{multline*}
\varphi(x) = \frac{-(\sigma x-1) +\Big( \sigma-1 -\frac{1}{2}(\sigma-1)\Big)\Big(\frac{\sigma x-1}{\sigma-1}\Big)^2 + O\Big((\sigma x-1)^3\Big)
}{
1-\sigma(1-x)
}\\
+ 
\frac{
\sigma x-1 -\frac{1}{2}(\sigma x-1)^2 +O\Big((\sigma x-1 )^3\Big)
}{
1-\sigma(1-x)
}\,.
\end{multline*}
This shows that the function~$\varphi$ and its derivative vanish at~$1/\sigma$ and 
\[
\varphi''\Big(\frac{1}{\sigma}\Big)= 2\sigma^2\frac{\sigma-1 -\frac{1}{2}(\sigma-1)-\frac{1}{2}(\sigma-1)^2}{(\sigma-1)^2(1-\sigma+1)}=
\sigma^2\frac{(\sigma-1)(2-\sigma)}{(\sigma-1)^2(2-\sigma)} = 
\frac{\sigma^2}{\sigma-1}.
\]
Therefore, we have 
\begin{equation}
\label{mvarphi}
m\varphi\Big((1-q)^\ell\Big) = m\Big((1-q)^\ell-1/\sigma\Big)^2\frac{\sigma^2}{2(\sigma-1)} + O\Big(m\big(\sigma (1-q)^\ell-1\big)^3\Big)\,.
\end{equation}
Writing this expression with the help of the variable~$\rho^*$ and replacing the first term with the estimates above, we obtain
\[
m\varphi\Big((1-q)^\ell\Big) =
m{\rho^*}^2 \frac{\sigma-1}{2}+O(m{\rho^*}^3)\,.
\]

The equivalence of the two times then leads to
\[
m{\rho^*}^2\, \frac{\sigma-1}{2} \mathop{\sim} \ell \ln \kappa\,,
\]
which means
\[
\Big(\sigma (1-q)^\ell -1\Big)^2= \frac{ \ell}{m} 2 (\sigma-1)\ln \kappa\, .
\]
We must therefore have 
\[
\sigma (1-q)^\ell= 1+\sqrt{\frac{ \ell}{m} 2(\sigma-1)\ln \kappa}\,.
\]
taking logarithm, we obtain thus
\[
q = 1-\exp\Bigg(\frac{-\ln \sigma}{\ell} +\frac{1}{\ell}\ln\bigg(1+\sqrt{\frac{ \ell}{m} 2 (\sigma-1)\ln \kappa}\bigg)\Bigg)\,.
\]
Expanding the exponential, we can develop the error threshold as 
\begin{multline*}
q^* = \frac{1}{\ell}\bigg(\ln \sigma - \sqrt{\frac{\ell}{m} 2 (\sigma-1)\ln \kappa} + \frac{\ell}{2m} 2 (\sigma-1)\ln \kappa - \cdots \bigg)\\
- \frac{1}{2\ell^2}\bigg(\ln \sigma - \sqrt{\frac{\ell}{m} 2 (\sigma-1)\ln \kappa} + \frac{\ell}{2m} 2 (\sigma-1)\ln \kappa - \cdots \bigg)^2
+ \cdots
\end{multline*}
The second term in the development after~$\ln \sigma /\ell$ must thus be of order~$1/\sqrt{\ell m}$. Suppose that we take~$q$ of the form
\[
q=\frac{\ln \sigma }{\ell}-\frac{c}{\sqrt{\ell m}}\,,
\]
for a certain positive constant~$c$. In this case, we have
\begin{align*}
(1-q)^\ell 
&= \frac{1}{\sigma} \exp\bigg( c\sqrt{\frac{\ell }{m}}+O\Big(\frac{1}{\ell}\Big)\bigg)\,.
\end{align*}
Since~$\ell^2$ is of order greater than~$m$, we get for~$m\big(\sigma(1-q)^\ell-1\big)^2$:
\begin{align*}
m\Big(\sigma(1-q)^\ell-1\Big)^2
&= c^2\ell +O\Big(\sqrt{\frac{\ell^3}{m}}\Big)\, .
\end{align*}
We must therefore compare the quantity~$c^2$ with 
$2 (\sigma-1)\ln \kappa.$
Since~$\ell$ dominates all the terms within the~$O$ in the expression of the persistence time of the previous theorem, we have proved theorem~\ref{seuilErreur}.
\subsection{In the case~$m/\ell \rightarrow \alpha$}
In this case, 
since the function~$\varphi$ is a bijection from~$[1/\sigma, 1]$ to~$[0, \ln\sigma]$, it is possible that the two times are of the same order but there is a condition on the limit~$\alpha$. We have that
\[
m\varphi\Big((1-q)^\ell)\Big) \mathop{\sim} \ell \ln \kappa\,,
\]
if and only if~$\alpha>\frac{\ln \kappa}{\ln\sigma}$.

And in this case we have that 
\[
(1-q)^\ell \rightarrow \varphi^{-1}\Big(\frac{\ln \kappa}{\alpha}\Big)\,.
\]
\subsection{In the case~$m/\ell \rightarrow 0$}
This time, no equivalence between the two times are possible, the neutral phase is always more stable than the quasispecies phase.

\section{In the case where~$\sigma e^{-a}-1=\sigma-2$}
\label{sigmaDeux}
First notice that hypothesis~\eqref{regimeconditionrho} leads to~$\sigma\geq 2$.
In this case, the quantity~$L_q$ tends towards~$0$, because the expession~\eqref{Lq} yields
\[
L_q = -\Big(\sigma (1-q)^\ell -1 -(\sigma-2)\Big)\,.
\]
Therefore, we cannot simplify the function~$F$ through~$\widetilde{F}$ and then bound~$\widetilde{F}$ uniformly, instead we will conduct the computations in a different way.
We thus define the new function~$F$ as
\begin{equation}
\label{F2}
F ( x ) = 
- \left( 1 - x \right) \ln \left( 1 - x \right) + x \ln \Big(\sigma(1-q)^\ell\Big) 
- \frac{\phi(x)\ln \phi(x)}{ L_q +(1-q)^\ell Q^\theta } \,.
\end{equation}
The critical point of function~$F$ defined by
\[
F'(x)=0\,,
\]
leads us to the equation
\begin{equation}
\label{eqPhi}
\sigma (1-q)^\ell (1-x) = \phi(x)\,,
\end{equation}
and to the expression
\[
\rho^* = 
\frac{\sigma(1-q)^\ell-1+(1-q)^\ell Q^\theta}{\sigma-1+(1-q)^\ell Q^\theta}\,.
\]
We develop this quantity in powers of~$Q^\theta$ and we write
\begin{equation}
\label{rhosigma2}
\rho^*= \rho^*_0
+\alpha Q^\theta \,,
\end{equation}
with
\[
\rho^*_0 = \frac{\sigma(1-q)^\ell-1}{\sigma-1}\,,
\]
and
\[
\alpha =
-\frac{\sigma(1-q)^\ell\big(1-(1-q)^\ell\big)}{(\sigma-1)^2}+O(Q^\theta) \,.
\]
We must then compute~$F(\rho^*)$, thanks to the identity~\eqref{eqPhi}, we have
\begin{multline}
\label{Frhosigma2}
F(\rho^*)=
\bigg(-1-\frac{1}{L_q+(1-q)^\ell Q^\theta}\bigg)
\ln(1-\rho^*)\\
-\bigg(\frac{1}{L_q+(1-q)^\ell Q^\theta}\bigg)
\ln\Big(\sigma(1-q)^\ell\Big)
+\frac{(1-q)^\ell Q^\theta\ln \phi(\rho^*)}{L_q+(1-q)^\ell Q^\theta}\,.
\end{multline}
Our goal is to recognize the function~$\varphi$ along with terms of order at most~$\frac{Q^\theta}{L_q(\sigma L_q+(\sigma-1)Q^\theta)}$.
Let us first handle the last logarithmic term 
\[
\frac{Q^\theta\ln \phi(\rho^*)}{L_q+(1-q)^\ell Q^\theta}
\bigg/
\frac{Q^\theta}{L_q(\sigma L_q+ (\sigma-1)Q^\theta)}
\mathop{\sim}\sigma L_q \ln\phi(\rho^*)\,.
\]
Since~$\phi(\rho^*)$ tends towards~$1$, this ratio is bounded.

We then develop the expression~\eqref{Frhosigma2} in powers of~$\rho^*$ and we get
\begin{multline*}
F(\rho^*)=
-\Big(1+\frac{1}{L_q}\Big)\ln\Big(1-\rho^*_0-\alpha Q^\theta\Big)
-\frac{1}{L_q}\ln\Big(\sigma (1-q)^\ell\Big)\\
+\Big(\frac{1}{L_q}-\frac{1}{L_q+(1-q)^\ell Q^\theta}\Big)
\ln\Big((1-\rho^*)\sigma(1-q)^\ell\Big)
+O\bigg(\frac{Q^\theta}{ L_q(\sigma L_q+(\sigma-1)Q^\theta)}\bigg)\,.
\end{multline*}
Since~$\rho^*$ tends towards~$(\sigma-2)/(\sigma-1)$, we have that 
\[
(1-\rho^*)\sigma(1-q)^\ell \rightarrow 1\,,
\]
the third term in the last expression is thus smaller than our remainder term. The expression of~$F(\rho^*)$ can then be written as 
\begin{equation*}
F(\rho^*) = \varphi\Big((1-q)^\ell\Big)
-\Big(1+\frac{1}{L_q}\Big)\ln\Big(1-\frac{\alpha{Q^\theta}}{1-\rho^*_0}\Big)
+O\bigg(\frac{Q^\theta}{L_q(\sigma L_q+(\sigma-1)Q^\theta)}\bigg)\,.
\end{equation*}
Finally, we have
\[
F(\rho^*) = \varphi\Big((1-q)^\ell\Big)
+O\bigg(\frac{Q^\theta}{ L_q(\sigma L_q+(\sigma-1)Q^\theta)}\bigg)\,,
\]
The sequel is quite analogous to the previous case.
\bibliographystyle{plain}
\bibliography{biblio}
\end{document}